\documentclass[12pt]{amsart}
\usepackage{epsfig,amssymb,epic,eepic,color}
\usepackage[all]{xy}
 \usepackage[T1]{fontenc}

\newcommand{\be}{\begin{enumerate}}
\newcommand{\ee}{\end{enumerate}}
\newcommand{\bi}{\begin{itemize}}
\newcommand{\ei}{\end{itemize}}

\def\R{\mathbb{R}}

\def\Z{\mathbb{Z}}

\def\Om{\Omega}

\def\ga{\gamma}

\def\be{\beta}
 
\def\de{\delta}

\def\vp{\varphi}

\def\si{\sigma}

\def\ep{\varepsilon}

\def\nd{\noindent}
\def\bull{\hfill$\Box$}

\begin{document}
\vskip 1cm
\begin{center}
{\sc A Morse complex on manifolds with boundary}
\end{center}

\title{}
\author{ Fran\c cois Laudenbach}
\address{Laboratoire de
math\'ematiques Jean Leray,  UMR 6629 du
CNRS, Facult\'e des Sciences et Techniques,
Universit\'e de Nantes, 2, rue de la
Houssini\`ere, F-44322 Nantes cedex 3,
France.}
\email{francois.laudenbach@univ-nantes.fr}
%\address{Universit\'e de Nantes, UMR 6629 du CNRS, 44322 Nantes, France}
%\email{francois.laudenbach@univ-nantes.fr}
%\maketitle



\thanks{supported by the French ANR program `Floer Power''}
%\date{March 18, 2010}

\begin{abstract} Given a compact smooth manifold $M$ 
with non-empty boundary and a Morse function, a pseudo-gradient
 Morse-Smale vector field adapted to the boundary
allows one to build a Morse complex whose homology is isomorphic to the 
(absolute or relative to the boundary) homology of $M$
 with integer coefficients. Our approach simplifies other methods which have 
been discussed in more specific geometric settings.

\end{abstract}
\maketitle

%\today
\thispagestyle{empty}
\vskip 1cm
\section{Introduction}
\medskip

We consider a smooth compact manifold $M$ of dimension $n$ with a non-empty 
boundary $\partial M$. 
A smooth function is said to be Morse when its critical points lie
in the interior of $M$, are not degenerate, and when its resriction
 to the boundary 
is Morse. The property of being Morse in this sense, with distinct
 critical values,
 is generic among the smooth functions on $M$. Let $f:M\to \R$ be
 a Morse function.
 There are two types of critical points of $f\vert\partial M$, 
 called type $N$  and D
 (we shall see later that N (resp. D) is like Neumann (resp. Dirichlet)):
 a critical point
 $p\in\partial M$
 is of type $N$ (resp. D) when $<df(p), n(p)>$ is negative (resp. positive),
 where $n(p)$ denotes the outward normal to the boundary at $p$.
 Looking at the change of  homotopy type of the sub-level set 
 $M^a:=\{x\in M\mid f(x)\leq a\}$ when $a$ is increasing, 
it is well known that 
change happens when $a$ crosses a critical value of $f\vert int(M)$ or
of $f\vert\partial M$ only when they are of type $N$; 
no change happens when crossing 
a critical value of type $D$.\\

 Thirty years ago, I tried to find Morse inequalities in this setting. 
Of course, $f$ gives 
 rise to a  Morse function $Df$ on the double manifold 
 $DM:= M\mathop{\cup}\limits_{\partial M}M$. But Morse inequalities for 
$Df$ are not sharp. At that time I did not succeed in finding
 a geometrical Morse complex 
 in the case of manifolds with boundary. I had even kept in mind
 the  idea that such a complex should  not exist. 
 Apparently, the problem seems to be still open, at least in the setting
 of generic assumptions on the boundary (see below).\\
 
 In a seminal paper \cite{witten}, E. Witten introduced 
a deformed Laplacian (now called the {\it  Witten Laplacian})
 on a closed Riemannian manifold equipped with
a Morse function and deduced in particular an analytic proof  
of the Morse inequalities. 
 Recently, Francis Nier explained  to me various works concerning the 
case of manifolds with non-empty boundary.  K. Chang and  J. Liu \cite{chang}
introduced two Witten complexes of differential forms 
constrained to satisfy some boundary conditions.
Two types of boundary conditions can be distinguished: Dirichlet conditions
 which cancel tangential components or Neumann conditions 
which cancel normal components.
% Indeed, there are two types of conditions: Dirichlet conditions 
%which are tangential, or Neumann conditions
% which are normal. 
These authors made the analysis easier
 by considering only flat metrics near the critical points.
These boundary problems with general metrics have been recently 
studied in \cite{nier} (Helffer-Nier) for the Dirichlet problem, 
in \cite{shubin} (Kolan-Prokhorenkov-Shubin) for Dirichlet and Neumann, and
in \cite{peutrec} (Le Peutrec) for the Neumann problem, with additional
 developments 
concerned with the asymptotic analysis.
%In  \cite{nier}  B. Helffer \& F. Nier   study the semiclassical limit 
%of the Witten Laplacian with Dirichlet boundary conditions 
%in the case of any Riemannian metric. 
%In  his recent thesis \cite{peutrec}, Dorian Le Peutrec solved the case of
%Neumann conditions.\\ 

 In case of the Dirichlet conditions \cite{chang} \cite{nier}, 
 when the deformation parameter $h$  is small enough,
  the De Rham complex $\Om^*$
 contains a finite dimensional sub-complex $F^*$ (of $\R$-vector spaces) 
 whose cohomology (when $M$ is orientable) is isomorphic 
to the relative cohomology
 $H^*(M,\partial M; \R)$. In that case, a basis of $F^*$ is in a bijective 
 correspondence with
 the critical points of $f$ in $int(M)$ and those of type $D$ 
on the boundary (up to a shift of their grading).
 In case of the Neumann boundary conditions \cite{chang} \cite{peutrec},
 a similar result holds:
  $F^*$ is generated  by the critical points of $f$ 
  in the interior and by those of type $N$ on the boundary;
  the cohomology of $F^*$ is isomorphic to the singular cohomology
 of $M$ with real coefficients.
  Hence, Morse inequalities follow. 
   Let us  also mention  work by M. Braverman and  V. Silantyev (\cite{braver})
 which is in the same spirit and deals with the Morse-Novikov theory.
But they introduce  some  extra 
condition which obliges them to exclude for 
instance the standard annulus in the plane equipped
 with the height function.

There are also approaches motivated by  Floer homology.
%\footnote{I thank 
%Jonathan Bloom who  informed  me of these works after my 
%first posted version.}.
 In \cite{akaho},
M. Akaho considers a very specific geometric situation along the boundary.
In \cite{kron}, P. Kronheimer and T. Mrowka consider functions and their
gradients having extensions to the double manifold $DM$ which are invariant 
by the canonical symmetry. 
Both settings %\cite{braver,akaho,kron} 
are  non-generic, and the setup and analysis are complicated 
by the presence of isolated trajectories in the boundary which preserves
Morse index.\\
   
These results (mainly  coming from analysis) beg the existence 
of a simpler geometric Morse complex associated to any Morse function 
on a manifold with boundary.
%it became ``necessary'' 
%to prove the existence of a geometric Morse complex in the setting
 %of manifolds with boundary.
In this note,  we are going to perform this program, working 
%in homology 
 with $\Z$-coefficients (or orientation-twisted coefficients)
 rather than with real cohomology as in Witten's work.
 We use the following notation:\bi
\item  $C_k$ denotes the set of critical points of
  $f:int(M)\to\R$ of index $k$;
  \item $N_k$ denotes the set of critical points 
of $f:\partial M\to \R$ of type $N$ and index $k$;
  \item $D_k$ denotes the set of critical points of 
$f:\partial M\to \R$  of type $D$ and index $k-1$ (notice that such
 a point has  index $k$ in the double manifold $DM$).
  \item $\vert\cdot\vert$ stands for the cardinality of the mentioned
 finite set.\\
 \ei
 
 \nd{\bf Theorem A.} {\it Let  $F^N_*$ be the free  graded $\Z$-module 
generated
 by $C_*\cup N_*$. 
  There exists  a differential $\partial:F^N_*\to F^N_{*-1}$,
 making $(F^N_*,\partial)$ a chain complex, such that
 the homology of $(F^N_*,\partial)$ is isomorphic to the singular
 homology
 $H_*(M,\Z)$.}\\
 
 As usual, Morse inequalities follow 
(see J. Milnor \cite{milnor}, p. 28, or R. Bott \cite{bott}, p. 338). 
They are contained
 in a polynomial identity. 
 Let $\mathcal M^N_f(T)$ be the  Morse polynomial of type $N$, where $T$ is 
the variable:
 $$\mathcal M^N_f(T)=\sum_k \vert C_k\cup N_k\vert \,T^k\,.$$
 Let $\mathcal P_M(T)$ be the Poincar\'e polynomial of $M$
 $$\mathcal P_M(T)=\sum_k {\rm rank}\,H_k(M;\Z)\, T^k\,.$$\\
 
 \nd {\bf Corollary A.}
 {\it We have $\mathcal M^N_f(T)-\mathcal P_M(T)=(1+T)Q^N(T)$,
 where $Q^N(T)$ is a polynomial with non-negative coefficients.}\\

 Making  now the critical points of type $D$ play the main role yields
 the following statement.\\
 
 \nd {\bf Theorem B.} {\it Let $F^D_*$ be the graded $\Z$-module generated by 
 $C_*\cup D_*$.
  There exists a co-differential $d:F^D_*\to F^D_{*+1} $,
 making $(F^D_*, d)$ a cochain complex, such that  
the cohomology of $(F^D_*, d)$ is isomorphic to the relative 
 cohomology $H^*(M, \partial M; \Z^{or})$ with coefficients twisted by the 
orientation
 of $M$.}\\
 
From this,  we deduce another family of Morse inequalities.  Let 
$\mathcal M^D_f(T)$ be the Morse polynomial
$$\mathcal M^D_f(T)=\sum_k \vert C_k\cup N_k\vert \, T^k.
$$
Let $\mathcal P_{(M,\partial M)}^{or}
$ be the relative Poincar\'e polynomial
$$\mathcal P_{(M,\partial M)}^{or}(T)
= \sum_k {\rm rank}\, H^k(M, \partial M; \Z^{or})\,T^k,$$
which is nothing but the symmetric polynomial of $\mathcal P_M(T)$:
$$ \mathcal P_{(M,\partial M)}^{or}(T)= T^n\, \mathcal P_M( 1/ T).$$\\

\nd {\bf Corollary B.} {\it We have 
$\mathcal M^D_f(T)-\mathcal P_{(M,\partial M)}^{or}(T)
=(1+T)Q^D(T)$, where $Q^D(T)$ is a polynomial with non-negative 
coefficients.}\\
 
Of course, corollaries A and B, which are  a direct consequence of 
theorems A and B respectively, are due  to Chang \& Liu \cite{chang}.\\

 I am very indebted to Francis Nier for encouraging me to return to this 
question and for his careful reading of a first version.
I also thank Jonathan Bloom and Claude Viterbo for useful information
and comments. \\
%I also thank Jonathan Bloom who informed me of other works in this topic.\\

 \section{Proof of Theorem A}

 \subsection{An adapted pseudo-gradient.}\label{adapted}
 We introduce a {\it pseudo-gradient} vector field $X$ for the Morse
 function $f$
 {\it adapted} to the boundary in the following sense (in case of 
closed manifolds it is a reformulation of \cite{meyer} by K. Meyer; 
a slightly more restrictive, but still generic, definition is given
 in \cite{slava}).
 \bi
 \item[1)] $X.f<0$ except at the critical points in $int(M)$ and
 at the critical points
 of type $N$ on the boundary;
 \item[2)] $X$ points inwards along the boundary except near the critical 
points
 of type $N$ where it is tangent to $\partial M$;
 \item[3)] if $p\in  int(M)$ is a critical point,  $X$ is  hyperbolic
 at $p$ and the quadratic form\break $\tilde q:=X^{lin}\,.\, d^2_pf $
is negative definite; here  $X^{lin}$ denotes the linear part of $X$ at $p$ 
(that is, $X^{lin}=DX(p)$) and  the second derivative $d^2_pf$ at $p$
 is thought of as 
a quadratic function defined near $p$ which is derived by $X^{lin}$;
 \item[4)] let $p\in\partial M$ be a critical point of type $N$; there are
 coordinates 
 $x=(y,z)\in \R^{n-1}\times\R$ of 
 $M$ near $p$, such that $M=\{z\geq 0\}$ and $f(x)=f(p)+z +q(y)$, where $q$ 
is a non-degenerate quadratic form; it is required that $X$ is a vector field
tangent to the boundary, vanishing %at $p$
 and hyperbolic at $p$, such that
$X^{lin}\,.\,\left(q(y)+z^2\right)$ is negative definite; 
 \item[5)] $X$ is Morse-Smale in the sense that the global 
unstable manifolds and
 the local
 stable manifolds are mutually transverse.\\
 \ei

\nd{\bf Proposition.} {\it Given the Morse function $f$, there exists an
adapted pseudo-gradient $X$.}\\

\nd {\bf Proof.} Look first at condition 4).
 The existence of coordinates where
$f$ reads as $f(x)=f(p)+z +q(y)$ mainly 
follows from  the Morse lemma with parameters (or
local stability of  Morse functions), up to the addition 
of a function depending on  $z$ only.  If $c(z)$ denotes the critical value of 
$y\mapsto f(y,z)$, the type $N$ assumption implies 
$ \frac{dc(z)}{dz}(0)>0$, hence the normal form holds after changing 
$z$ by applying 
the inverse function theorem. Moreover, an easy calculation shows
that, away from $p$,  $X^{lin}.\bigl(q(y)+z^2)\bigr)<0$  implies 
$X^{lin}.\bigl(q(y)+z)\bigr)<0$. 
Therefore, $X$ is a pseudo-gradient for $f$, that is
 $X.f<0$ near $p$ (except at $p$).
 
Now  the local existence of $X$ is clear 
 and local pseudo-gradients can be glued together by a partition of unity.
 By construction, $X$ is positively complete. 
Let $X^t$ denote 
its flow; it is defined
 on an open set of $M\times \R$ containing 
 $M\times[0,+\infty)$. The global unstable manifold $W^u(p)$ of a critical 
point $p$ of index 
$k$ is the image of a non-proper embedding of $\R^k$ into $M$.

There is also 
 a local stable manifold $W^s_{loc}(p)$ which is diffeomophic to $\R^{n-k}$
 when 
$p\in int(M)$, or  to $\R^{n-k}\cap\{z\geq 0\}$ when $p$ is a type $N$ 
critical point in
the boundary (here  $\R^{n-k}$ is a
 space on which $q(y)+z^2$ is positive definite).
 It is properly embedded when it is 
 truncated to the sub-level set $M^{f(p)+\ep}$. 
Following Smale \cite{smale}, condition 5) is generically
fulfilled among the vector fields meeting conditions 1)-4). \bull\\

If $p$ has index $k$ the frontier 
of $W^u(p)$, that is the set of points in the closure which are not in 
$W^u(p)$, is contained in the union of the unstable manifolds of critical 
points of index less than $k$. An orientation of 
$W^u(p)$ is  chosen arbitrarily. Then  $W^s_{loc}(p)$ 
is co-oriented by the orientation of $W^u(p)$.\\

Given a pair of critical points $(p,q)$ respectively of index $k$
and $k-1$, there are only finitely many flow lines (up to translation in time)
$X^t(x), t\in\R$, such that $X^t(x)\to p$ as $t\to -\infty$ and 
$X^t(x)\to q $ as 
$t\to +\infty$.  Each such line has a sign according to the co-orientation of
$W_{loc}^s(q)$ with respect to the orientation of $W^u(p)$.
Let $m_{pq}$ be the algebraic sum of these signs, summing up over all 
connecting orbits from $p$ to $q$.\\

\subsection{The boundary morphism.} \label{morphism}We are going to define 
$\partial: F_k\to F_{k-1}$ {\it \`a la Witten-Floer} (see \cite{witten,floer}).
Let $p\in C_k\cup N_k$ be a generator. The boundary morphism is defined by:
$$\partial <p>=\sum_q m_{pq}<q>,$$
where $q$ runs in $C_{k-1}\cup N_{k-1}$.\\

\nd {\bf Proposition.} {\it  We have $\partial\circ\partial=0$.}\\

\nd {\bf Proof.} Consider  $p$, a critical point of index $k$, and $q$,
 a critical point 
of index $k-2$. We have to prove that $<\partial\circ\partial (p), q>=0$.
 In \cite{laudenbach} it is given a complete description of
 the closure of 
$W^u(p)$; only the effect of critical points of index $k-1$ is useful for
 yielding the following description. In a level set  $L:= \{f(x)=f(q)+\ep\}$,
 we look at the trace $A$ of
$W_{loc}^s(q)$, a $(n-k+1)$-dimensional sphere or   proper disk,
and the trace $B$ of the closure of $W^u(p)$. The intersection $A\cap B$ 
lies in 
the interior of $L$ since $W^u(p)$ lies in $int(M)$ except near $p$ 
when $p\in\partial M$.
 It is made of simple closed curves and  arcs having end points in common
  and whose interior are mutually disjoint.
The closed curves and 
 the open arcs correspond to connecting orbits from $p$ to $q$; 
the end points correspond to broken connecting orbits going through some 
critical points of index $k-1$.  Notice that the sign of $m_{pz}m_{zq}$ 
does not depend on the chosen  
orientation of the unstable manifold of $z$, 
an index $(k-1)$ point connected to $p$
and to $q$.
Each open arc in $A\cap B$ is oriented.
 Then one of his end points is equipped with $+$,
the other with $-$, and the sum of all these signs 
is $<\partial\circ \partial (p), q>$. Hence, it is zero.\bull\\

\subsection{Invariance.} \label{invariance} 
(This property is not needed in the proof of theorems A and B,
 but it is interesting in itself). 
Of course, the complex $(F^N,\partial)$ depends 
on the Morse function $f$ and on an adapted pseudo-gradient $X$.
Given two pairs $(f_0, X_0)$ and $(f_1, X_1)$ where, for $i=0, 1$,
 $f_i$ is a Morse function and $X_i$ is an adapted pseudo-gradient,
 they are connected
 by a path of pairs $(f_t, X_t),\ t\in [0,1]$, if we allow us to 
cross ``codimension 1 accidents'' which arrive at
 finitely many times $t_1, t_2,...t_\ell$. When $t$ is distinct of 
those times, $(f_t, X_t)$ is a pair of a Morse function and an adapted 
 pseudo-gradient.
 The list of the possible accidents is yielded by appropriate
 transversality theorems. We list below the ``codimension 1 accidents'' 
for a  function 
 $f$;
 and those of pseudo-gradient $X$. 
 Here, it is convenient to set $f_\partial:=f\vert \partial M$.
 
 \bi
 \item[f1)] $f$ has a degenerate critical point at $p\in int(M)$ and
 rank\,$d^2f(p)=n-1$;
 \item[f2)] $f_\partial$ has a degenerate critical point at $p\in \partial M$,
 rank\,$d^2f_\partial(p)= n-2$ and $df(p)\not=0$;
  \item[f3)] $df(p)$ vanishes at some point $p\in \partial M$,
  $d^2f(p)$ and $d^2f_\partial(p)$ are both non-degenerate and 
 have the same index; 
  \item[f4)] $df(p)$ vanishes at some point $p\in \partial M$,
  $d^2f(p)$ and $d^2f_\partial(p)$ are both non-degenerate and 
  index\,$d^2f(p)$= index\,$d^2f_\partial(p)+1$.
  \ei
  When crossing an accident f1), a pair of critical points of consecutive 
indices 
  is created/cancelled %for $f$ 
in  $int(M)$. When crossing an accident f2),
  a pair of critical points of consecutive indices and of the same 
type $N$ or $D$
  is created/cancelled for $f_\partial$ in  $\partial M$.
  Crossing an accident f3) can be modelled as follows: there are local
 coordinates
  near $p$, $x=(y,z)\in \R^{n-1}\times \R$  where $M=\{z\geq 0\}$ and
  $f(x)= f(p)+z^2+q(y)$ where $q$ is a non-degenerate quadratic form 
(say of index $k$). 
  For crossing the accident one leaves this function fixed and translates 
$M$ by
  $(y,z)\mapsto (y, z+t)$, $t\in[-\ep,+\ep]$. For $t=\ep$, we have a point 
  of type $N$ and index $k$ on the boundary; for $t=-\ep$, we have a point of
 type $D$ 
  on the boundary and a new point of index $k$ in $int(M)$.
  The model for crossing f4) is similar but  $f(x)= f(p)-z^2+q(y)$. 
When moving from
  $t=\ep$ to $t=-\ep$, a point of type $D$ and index $k$ on the boundary is 
replaced
  by a point of type $N$ and index $k$ on the boundary and a point of 
index $k+1$
  in the interior.\\
  
  We now list the ``codimension 1 accidents'' for a pseudo-gradient $X$; 
they never
  happen at the same time as a ``codimension 1 accident'' of the function
 whose $X$
  is a pseudo-gradient.
  \bi
  \item[g1)] creating/cancelling a pair of connecting orbits from a point 
of index k %(in $int(M)$ or $\partial M$)
 to a point of index $k-1$ (the considered critical points may lie in $int(M)$
or in $\partial M)$;
  \item[g2)] there is an orbit connecting two points of the same index.\\
  \ei
  
 \nd {\bf Proposition.} {\it Let $(F_*^i,\partial ^i)$ be the complex 
associated to the pair 
 of a Morse function and a pseudo-gradient
 $(f_i, X_i)$, $i=0,1$. Then they are quasi-isomorphic: there is a chain 
morphism
 from one to the other inducing an isomorphism in homology.}\\
 
 \nd {\bf Proof. }
According to the preceding discussion, we may assume that there exists a path
 $(f_t, X_t)$ with one accident only. The proposition is proved
  by examining  each accident.
  When crossing f1), an acyclic complex of rank 2 
(that is, $0\to\Z\mathop{\longrightarrow}\limits^{\cong}\Z\to 0$)
 is added or removed.
 The same happens
  when crossing f2), if the considered points are of type $N$; if they are of 
type $D$ the complex is unchanged. When crossing f3), the complex 
remains  unchanged. 
  When crossing f4), again an acyclic complex of rank 2 is added or removed.
  Concerning the accidents of $X$, there are those which are encountered in the
  analogous discussion for closed manifolds (see \cite{laudenbach}). In each
 case we get the desired quasi-isomorphism. \bull\\

 \subsection{The homology of $(F^N_*,\partial)$}\label{homology}
 
 ${}$
 \medskip

 \nd {\bf Proposition.} {\it The homology of $(F^N_*,\partial)$ is isomorphic 
to 
 $H_*(M; \Z)$.}\\
 
 \nd {\bf Proof.} %After the invariance proposition proved in \ref{invariance}, it is sufficient to consider a particular Morse function. 
According to the next lemma, it is 
allowed to assume that there are no critical points of type $N$.
 Let us achieve 
the proof under this assumption. In that case the complex does not ``see'' the
 boundary  since
 all  connecting orbits lie in $int(M)$. 
Even, the global unstable manifolds lie
 in the interior of $M$.
 So, we can deal with such a function as on a closed manifold. For instance, 
 an adapted pseudo-gradient being chosen, it is possible to reorder the
 critical values
 so that the function becomes self-indexing (the critical value of a critical
 point is its index). Indeed, we recall the following fact in Morse theory
 (see Cerf \cite{cerf}, II 2.3). 
\bi
\item[] {\it Let $(f,X)$ be a Morse function and an adapted pseudo-gradient.
 Let $p$ 
and $q$ two critical points with $f(p)>f(q)$. Assume that the open interval
$\bigl(f(q),f(p)\bigr)$ contains no critical value and that there are no  
connecting orbits
from $p$ to $q$. Then there exists a path of Morse functions $f_t,\ t\in[0,1]$,
with $f_0=f$, $f_1(p)<f_1(q)=f(q) $ and $X$ is a pseudo-gradient for every 
$f_t$.
}
\ei
Once $f$ is self-indexing, $M$ has the homotopy type of a CW-complex, 
with one cell of dimension $k$ for each critical point of index $k$
(see Milnor \cite{milnor3}, theorem 7.4, p. 90). In that case, our Morse chain 
complex is exactly the so-called {\it cellular chain complex}. 
By the {\it cellular homology theorem} (see Milnor \cite{milnor2}, theorem A.4 
p. 263), its homology is the  
homology of $M$. So the proof will be finished after the following lemma.\\

\nd {\bf Lemma.} {\it Let $(f,X)$ be a pair of a Morse function on $M$ and
an adapted pseudo-gradient. Let  $p$ be a critical point 
of type $N$  in $\partial M$ and $U$ be a neighborhood of $p$.
Then there exist a  path of functions
$f_t$ and an isotopy $\vp_t$ of embeddings $M\to M$,
 $t\in[0,1],$ with $f_0=f, \ \vp_0=Id$, 
satisfying the following properties :  
\bi
\item[1)] the support of both $f_t-f$ and $\vp_t $ is contained in $U$; 
 \item[2)] $f_t$ is Morse when $t\not= 1/2$; 
 \item[3)] $f_1$ has one critical point $p_1$  in $U\cap int( M)$;
\item[4)] $f_1\vert \partial M$ has   one critical point in $U
\cap\partial M$ and it is of type $D$;
\item[5)] $\vp_{1}\left(W^u(p)\right)$ is the unstable manifold of
 $p_1$ for an adapted 
pseudo-gradient $X_1$ of $f_1$, whose other invariant manifolds are those
 of $X$.
In particular, the  Morse complexes associated to $X$ and $X_1$ are the same.\\
\ei}

Of course the accident at time $t=1/2$ is f3) in list \ref{invariance}.\\
 
 \nd {\bf Proof.} We start with a model for $f$  in coordinates 
 $x=(y,z)\in \R^{n-1}\times \R$ on a small open neighborhood $U$ of 
 the critical 
point $p$:
 $f(x)= f(p)+ z+ q(y)$, where $q(y)$
is a non-degenerate quadratic form and $M\cap U\subset \{z\geq 0\}$.
We consider $f_{1/2}$ defined on $U$ by $f_{1/2}= f(p)+q(y)+z^2$. 
 Let $\rho: [0,\de]\to [0,1]$ be a smooth bump function: 
$\rho(s)=1$ for $s$ close to 0
and $\rho(s)=0$ for $s$ close to $\de$; here $\de$ is small enough so that
 the support
of $(y, z)\mapsto \bigl(\rho(\Vert y\Vert), \rho(z)\bigr)$ is a compact 
set in $U$. 
Without loss of generality, we may assume that the support of $\rho$ 
contains in its interior the connected domain of $\partial M$ 
along which  $X$ is tangent to the boundary.

 For $\ep>0$ small enough, we look at the 
restriction 
of $f_{1/2}$ to $\{z\geq \ep\rho(\Vert y\Vert)\}$. It is conjugate
 to  $f\vert M\cap U$
 because 
both functions have no critical points  in $U$ and have restrictions to
 $\partial M\cap U$ which are Morse with one critical point of the same index.
 Let $M_0$ be the closure of
 $M\setminus \{0\leq z\leq \ep\rho(\Vert y\Vert)\}$ 
 and $\psi: M_0\to M$ be 
 a diffeomorphism such that both functions  $f_0:=\psi^*f$ and $f_{1/2}$ 
 have the same germ along $\{z=\ep\rho(\Vert y\Vert)\}$. 
 For $t\in [0,1]$, introduce the 
 manifold $M_t$ and its function $f_t$ obtained from $(M_0, f_0)$ by gluing 
 $\{\ep(1-2t)\rho(\Vert y\Vert)\leq z\leq \ep\rho(\Vert y\Vert)\}$
 endowed with $f_{1/2}$.
 Up to a diffeomorphism $\psi_t :M_t\to M$, this path $f_t$ is the desired one.\\
 
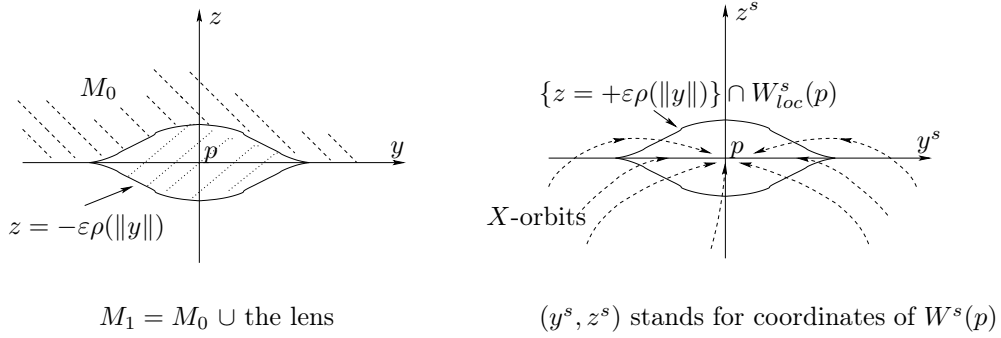
 \begin{center}
${}\quad$\begin{picture}(0,0)%
\epsfig{file=figure1.pstex}%
\end{picture}%
\setlength{\unitlength}{1579sp}%
\begingroup\makeatletter\ifx\SetFigFont\undefined%
\gdef\SetFigFont#1#2#3#4#5{%
  \reset@font\fontsize{#1}{#2pt}%
  \fontfamily{#3}\fontseries{#4}\fontshape{#5}%
  \selectfont}%
\fi\endgroup%
\begin{picture}(6538,5089)(2326,-6548)
\put(3751,-6436){\makebox(0,0)[lb]{\smash{{\SetFigFont{10}{12.0}{\familydefault}{\mddefault}{\updefault}{\color[rgb]{0,0,0}$M_1=M_0\ \cup\ $the lens}%
}}}}
\put(5476,-1711){\makebox(0,0)[lb]{\smash{{\SetFigFont{10}{12.0}{\familydefault}{\mddefault}{\updefault}{\color[rgb]{0,0,0}$z$}%
}}}}
\put(8326,-3736){\makebox(0,0)[lb]{\smash{{\SetFigFont{10}{12.0}{\familydefault}{\mddefault}{\updefault}{\color[rgb]{0,0,0}$y$}%
}}}}
\put(2326,-5011){\makebox(0,0)[lb]{\smash{{\SetFigFont{10}{12.0}{\familydefault}{\mddefault}{\updefault}{\color[rgb]{0,0,0}$z=-\ep\rho(\Vert y\Vert)$}%
}}}}
\put(5401,-3811){\makebox(0,0)[lb]{\smash{{\SetFigFont{10}{12.0}{\familydefault}{\mddefault}{\updefault}{\color[rgb]{0,0,0}$p$}%
}}}}
\put(3451,-2836){\makebox(0,0)[lb]{\smash{{\SetFigFont{10}{12.0}{\familydefault}{\mddefault}{\updefault}{\color[rgb]{0,0,0}$M_0$}%
}}}}
\end{picture}%
${}\quad\quad$\begin{picture}(0,0)%
\includegraphics{figure2.pstex}%
\end{picture}%
\setlength{\unitlength}{1579sp}%
\begingroup\makeatletter\ifx\SetFigFont\undefined%
\gdef\SetFigFont#1#2#3#4#5{%
  \reset@font\fontsize{#1}{#2pt}%
  \fontfamily{#3}\fontseries{#4}\fontshape{#5}%
  \selectfont}%
\fi\endgroup%
\begin{picture}(8055,5160)(1576,-6619)
\put(5401,-3811){\makebox(0,0)[lb]{\smash{{\SetFigFont{10}{12.0}{\familydefault}{\mddefault}{\updefault}{\color[rgb]{0,0,0}$p$}%
}}}}
\put(8326,-3736){\makebox(0,0)[lb]{\smash{{\SetFigFont{10}{12.0}{\familydefault}{\mddefault}{\updefault}{\color[rgb]{0,0,0}$y^ s$}%
}}}}
\put(1576,-4936){\makebox(0,0)[lb]{\smash{{\SetFigFont{10}{12.0}{\familydefault}{\mddefault}{\updefault}{\color[rgb]{0,0,0}$X$-orbits}%
}}}}
\put(2401,-2986){\makebox(0,0)[lb]{\smash{{\SetFigFont{10}{12.0}{\familydefault}{\mddefault}{\updefault}{\color[rgb]{0,0,0}$\{z=+\ep\rho(\Vert y\Vert)\}\cap W^ s_{loc}(p)$}%
}}}}
\put(5476,-1711){\makebox(0,0)[lb]{\smash{{\SetFigFont{10}{12.0}{\familydefault}{\mddefault}{\updefault}{\color[rgb]{0,0,0}$z^s$}%
}}}}
\put(2401,-6511){\makebox(0,0)[lb]{\smash{{\SetFigFont{10}{12.0}{\familydefault}{\mddefault}{\updefault}{\color[rgb]{0,0,0}$(y^s,z^ s)$ stands for coordinates of $W^ s(p)$}%
}}}}
\end{picture}%

\end{center}

\medskip

 More precisely, the isotopy $\psi_t, \ t\in[0,1],$
 can be chosen so that $\psi_{1/2}= Id$.
 %and $\psi_t^*f=f_t$.
  Recall that $X$ is a pseudo-gradient  in $U$ for both functions $f$ and 
   $f_{1/2}$ (see \ref{adapted}). At $t= 1/2$, the deformation of functions
 is stopped
   and a small isotopy $\chi_t$, is applied for pushing
 $W^u(p)\setminus\{p\}$ into $int(M)$; it is chosen  to be
supported in  $int(M_1)$ and to  leave $W_{loc}^s(p)$ fixed.
  As this isotopy is $C^\infty$ small,   $X_1:=\chi_{1*}X$ is still
 a pseudo-gradient of $f_{1/2}$ and 
 $\chi_1(W^u(p))$ is its  unstable
 manifold of $p$. 

  Moreover, $X_1$ points 
  toward  $int(M)$ along $U\cap\{z=0\}$ except along  $W_{loc}^s(p)$.
  If $\eta>0$ is small, $X_1$ points inwards 
along $\partial M_{1/2+\eta}\cap U$.  
  Therefore, if the path of deformation $(M_t,f_t)$ 
 is stopped at time $t=1/2+\eta$
  we already get the desired conclusion.
 \bull\\

\nd {\bf Remark.} (C. Viterbo) Instead of deforming $(f,X)$ as in the 
preceding lemma, one could invoke Conley's theory of {\it isolating blocks}
\cite{conley} which in our setting states as follows (compare \cite{floer} section 3 for a similar situation, but our flow is not negatively complete).

%\bi
%\item[] 
{\it There exists a finite 
filtration
$M_{-1}=\emptyset\subset M_0\subset \ldots\subset M_k\ldots \subset M_n=M$
which is positively invariant by the flow $X^t$ and
such that $M_k\setminus M_{k-1}$ contains all critical points in $C_k\cup N_k$
and every positive orbit in $M_k$  which does not go to a critical point
enters $M_{k-1}$ without getting out of.}\\

  Our $F_k$ can be identified to $H_k(M_k, M_{k-1};\Z)$ and our 
$\partial$ is nothing but the connecting morphism 
$H_k(M_k, M_{k-1};\Z)\to H_{k-1}(M_{k-1}, M_{k-2};\Z) $ in the long exact
 homology sequence of the triple $(M_k, M_{k-1},M_{k-2})$ (as in Milnor 
\cite{milnor2} A.4).

\subsection{Complement}\label{complement}

We have statements similar to theorem A for every local system of coefficients
(or flat bundle with discrete fiber). Such a bundle is given by a
 representation
of the fundamental group $\pi:=\pi_1(M,x_0)$ (where $x_0$ is a base point) into
an Abelian discrete group. For instance, $\Z^{or}$ is defined 
by the first Stiefel-Whitney
class of $\tau M$, the tangent bundle of $M$: $w_1:\pi\to \Z/2$.
The complex $F_*^{N,or}$ has the same graded module. But the differential
is changed. At the beginning some arc is chosen from the base point to each
 critical point. Therefore a connecting orbit $\ga$ from $p$ to $q$ defines a
 loop $\tilde\ga$
and the contribution of $\ga$ to $m_{pq}$ is twisted by $w_1(\tilde\ga)$. So
 we have a new differential $\partial^{or}$. The proof that 
$\partial^{or}\circ\partial^{or}=0$ is alike.
In that case, theorem A states as the following: {\it The homology
of $(F_*^{N,or}, \partial^{or})$ is isomorphic to $H_*(M;\Z^{or})$.}
New Morse inequalities follow.  A good example is given by the compact 
M\"obius band
in $\R^3$ with a height function. In that case the Betti numbers of
$H_*(M;\Z^{or})$ are: $b_0=0,b_1=0, b_2=0$; in relationship to theorem B,
we also list the Betti numbers of the relative homology 
$H_*(M,\partial M;\Z^{or})$: $b_0= 0, b_1=1, b_2=1$.   \\

\section{Proof of Theorem B and complements}

\subsection{Upside down.}
We look at the Morse theory of the Morse function $-f$. The critical points are
 the same
but their indices and types are changed.  A critical point in $int(M)$ of index
 $k$ for 
$f$ has index $n-k$ for $-f$. A critical point in $\partial M$ which is of type
 D and index 
$k-1$ for $f_\partial:= f\vert \partial M$ 
(then it belongs to $D_k$) is of type $N$ and index $n-k$ 
for $-f_\partial$. A critical point in $\partial M$ which is of type $N$ and
 index 
$k$ for $f_\partial$ is of type $D$ and index $n-k-1$  
for $-f_\partial$; then it belongs to $N_{n-k}(-f)$.

There is another change to make. Indeed, the vector field $-X$ is not adapted 
as it points outwards along the boundary except near the critical points of
 former type $N$.
So we appeal another vector field $X^-$ which is adapted to $-f$.
With the data $(-f, X^-) $ we can form a Morse complex. 
But we insist to keep the initial grading, that is: 
$F^D_k$ is freely generated by $C_k\cup D_k$, a set of critical points of 
index $n-k$
for $-f$. The pseudo-gradient  $X^-$ yields a differential which is of
 degree +1 in this grading: it is a co-differential and the complex  is a 
cochain complex. According to 
theorem A its (co-)homology is isomorphic to $H_*(M,\Z)$, up to the grading which is changed into the complementary one. Therefore, it is the Poincar\'e dual
of the homology, hence the cohomology with integer 
coefficients twisted by the orientation.
Thus theorem B is proved.
\bull\\

As we did in a variant of theorem A in 
\ref{complement}, we can twist the  codifferential $d$ by the orientation. 
We get a cochain complex  $(F_*^{D,or}, d^{or})$
whose cohomology is isomorphic to $H^*(M,\partial M; \Z)$.
From this we get another Morse inequalities involving the same Morse 
polynomial:
$$\mathcal M^D_f(T)-\mathcal P_{(M,\partial M)}(T)=(1+T)Q'(T).
$$
We recall here, that homology and cohomology have the same Poincar\'e
 polynomial,
as the free factor of $H^k(.)$ is the dual of the free factor of $H_k(.)$, 
as long as the same system of local coefficients is used for both.

\subsection{Morse inequalities from the double.}
On the one hand, the set $C_k(Df)$
of  critical points of index $k$ for the ``double'' Morse function $Df$ defined
on the double $DM$ of $M$  consists of two copies of $C_k$, one copy of $N_k$
and one copy of $D_k$.
On the other hand, $H_k(DM; \Z)\cong H_k(M;\Z)\oplus H_k(M,\partial M; \Z)$.
Indeed, the inclusion $M\hookrightarrow DM$ gives rise to 
a long exact sequence
$$\ldots \to H_k(M)\to H_k(DM)\to H_k(M,\partial M)\to H_{k-1}(M)\to \ldots\,;
$$
but this sequence splits since any relative cycle, glued to its symmetric one, 
 gives a cycle in $DM$. Thus the Morse inequalities for $Df$ is just obtained
 by adding
both  Morse inequalities obtained 
previously, the one for $\mathcal M^N_f -\mathcal P_M$ and the one for 
$\mathcal M^D_f-\mathcal P_{(M,\partial M)}$:
$$
\mathcal M_{Df}(T)- \mathcal P_{DM}(T)=(1+T)Q''(T).
$$
Thus, it is weaker than its summands. 

\subsection{Relative homology and its pairing with absolute homology}
\footnote{\ I thank M. Akaho for 
 asking me about the relative homology.}
Here we need the following notation: for $p\in C_k\cup D_k$, $W^u(p,X^-)$
will 
denote the unstable manifold of the pseudo-gradient vector field 
$X^-$ adapted to the function $-f$; it is equipped with an orientation.
If we want to think of it as a relative $k$-cochain, it is necessary to endow
$p $ with  a {\it local orientation } of $M$. Then, if $\si$ is an oriented 
$k$-simplex tranverse to  $W^u(p,X^-)$ it is allowed to count 
the intersection points in  $\si\cap W^u(p,X^-)$ with signs.
It is the way of thinking of the complex $(F_*^D,d)$ as a cochain complex
with orientation-twisted coefficients calculating $H^*(M,\partial M;\Z^{or})$.

Now, we consider the dual complex $(\check{F}_*^D, \check{d})$ in which 
the differential is the transpose of $d$. Then its homology 
is isomorphic to $H_*(M,\partial M;\Z^{or})$; here we use the 
universal-coefficient theorem for cohomology (\cite{spanier}, theorem 5.5.3)
and the finite generation of $(F_*^D,d)$.
 Geometrically, 
$\check{d}$ can be described as follows. Let $p\in C_k\cup D_k$ and 
$q\in C_{k-1}\cup D_{k-1}$. We have to calculate $\check{m}_{pq}:=
<\check{d}(p),q>$. It is the algebraic number of connecting orbits of $X^-$
from $q$ to $p$ counted with the orientation of $W^u(q,X^-)$
and the co-orientation of $W^s_{loc}(p,X^-)$.\\

For making the pairing $H_k(M,\partial M;\Z^{or})\otimes 
H_{n-k}(M; \Z)\to\R$ explicit, we need to ask a general position of $X^-$
with respect to $X$: for each pair of integers $(k,k')$
and for each $p\in C_k\cup D_k$ and 
$p'\in C_{n-k'}\cup N_{n-k'}$ the unstable manifold $W^u(p,X^-)$ is transverse
to the unstable manifold  $W^u(p',X)$. When $k=k'$, it is easy to see that
only compact parts of these invariant manifolds could intersect; therefore,
there are finitely many intersection points. The sign is calculated by using 
three data; the orientation of $W^u(p',X)$, the orientation of
$W^u(p,X^-)$ and the local orientation of $M$ at $p$. 

Arguing  this  way, we avoid to identify  $\check{F}_*^D$ to a quotient
complex of $F_*^N$ by a subcomplex isomorphic to the Morse complex 
of $f\vert\partial M$. Such an identification should require difficult formulas
similar to those in \cite{kron}, section 2.4.

\end{document}